\documentclass[runningheads,a4paper]{llncs}

\usepackage{amssymb}
\usepackage{graphicx}

\usepackage{amsmath}
\usepackage{algpseudocode}
\usepackage{algorithm}
\usepackage{url}
\urldef{\mailsa}\path|{d.andreagiovanni,mett,pulaj}@zib.de|
\newcommand{\keywords}[1]{\par\addvspace\baselineskip
\noindent\keywordname\enspace\ignorespaces#1}

\begin{document}

\mainmatter

\title{An (MI)LP-based Primal Heuristic for 3-Architecture Connected Facility Location\\in Urban Access Network Design\thanks{
This is the authors' final version of the paper published in: Squillero G., Burelli P. (eds), EvoApplications 2016: Applications of Evolutionary Computation, LNCS 9597, pp. 283-298, 2016.
DOI:  10.1007/978-3-319-31204-0\_19.
The final publication is available at Springer via http://dx.doi.org/10.1007/978-3-319-31204-0\_19 \; 
The work of Fabio D'Andreagiovanni and Jonad Pulaj was partially supported by the \emph{Einstein Center for Mathematics Berlin} (ECMath) through Project MI4 (ROUAN) and by the \emph{German Federal Ministry of Education and Research} (BMBF) through Project VINO (Grant 05M13ZAC) and Project \emph{ROBUKOM} (Grant 05M10ZAA)\cite{BaEtAl14}.}}

\titlerunning{An (MI)LP-based heuristic for 3-Architecture Connected Facility Location}

\author{Fabio D'Andreagiovanni$^{1,2}$
\and
Fabian Mett$^{1}$
\and
Jonad Pulaj$^{1}$}

\authorrunning{F. D'Andreagiovanni, F. Mett, J. Pulaj}

\institute{
$^{1}$Dept. of Mathematical Optimization, Zuse Institute Berlin (ZIB), Takustr. 7, 14195 Berlin, Germany\\
%$^{2}$DFG Research Center MATHEON, Technical University Berlin, Stra{\ss}e des 17 Juni, 10623 Berlin, Germany\\
%$^{3}$Einstein Center for Mathematics (ECMath), Stra{\ss}e des 17 Juni, 10623 Berlin, Germany\\
$^{2}$Institute for System Analysis and Computer Science, National Research Council of Italy (IASI-CNR), via dei Taurini 19, 00185 Roma, Italy\\
\mailsa\\
}

\toctitle{Lecture Notes in Computer Science}
\tocauthor{Authors' Instructions}
\maketitle

\begin{abstract}
We investigate the 3-architecture Connected Facility Location Problem arising in the design of urban telecommunication access networks integrating wired and wireless technologies. We propose an original optimization model for the problem  that includes additional variables and constraints to take into account wireless signal coverage represented through signal-to-interference ratios.
Since the problem can prove very challenging even for modern state-of-the art optimization solvers, we propose to solve it by an original primal heuristic that combines a probabilistic fixing procedure, guided by peculiar Linear Programming relaxations, with an exact MIP heuristic, based on a very large neighborhood search. Computational experiments on a set of realistic instances show that our heuristic can find solutions associated with much lower optimality gaps than a state-of-the-art solver.
\keywords{Telecommunications, FTTX Access Networks, Connected Facility Location, Mixed Integer Linear Programming, Tight Linear Relaxations, MIP Heuristics.}
\end{abstract}

\section{Introduction}

In the last two decades, telecommunications have increasingly assumed a major role in our everyday life and the volume of traffic exchanged over wired and wireless networks has enormously increased.
%, causing a growth of networks in size and complexity.
Major telecommunications companies forecast that such growth
%, also induced by the continuously growing need for faster and more sophisticated telecommunication services,
will powerfully continue, thus requiring the need for more technologically complex networks.
In this context, \emph{telecommunication access networks}, which connects users to service providers, have become a vital part of urban metropolitan infrastructures. A critical component of such networks is represented by optical fiber connections, which provide higher capacity and better transmission rates than the old copper-based connections. In the last years, the trend in access networks has been to provide broadband internet access through different types of optical fiber deployments. These several deployments, usually called \emph{architectures}, are denoted as a whole by the acronym \emph{FTTX} (Fiber-To-The-X), where the \emph{X} is  specified on the basis of where the optical fiber granting access is terminated: major examples of architectures are \emph{Fiber-To-The-Home (FTTH)}, bringing a fiber directly to the final user, and \emph{Fiber-To-The-Cabinet (FFTC)} and \emph{Fiber-To-The-Building} (FTTB), respectively bringing a fiber to a street cabinet or to the building of the user and then typically connecting the fiber termination point to the user through a copper-based connection. We refer the reader to \cite{GrRaWe14} for an exhaustive introduction to FTTX network design and to \cite{ZoMyTo15} for a thorough discussion about the features of FTTH network design. Nowadays, an access network implementing a \emph{full} FTTH architecture seems impractical, because of its extremely high deployment costs and since not all users are willing to pay higher fees for faster connections. As a consequence, in recent times higher attention has been given to deployments mixing \emph{two architectures} like FTTH and FTTC/FTTB (e.g., \cite{LeEtAl13}). An even more recent and promising trend has been represented by the integration of wired and wireless connections, providing service to users also through wireless links and leading to \emph{3-architecture} networks that includes also the so-called \emph{Fiber-To-The-Air} (FTTA) architecture \cite{GhMaAs09,GrRaWe14}. This integration aims to get the best of both worlds: the high capacity offered by optical fiber networks and the mobility and ubiquity offered by wireless networks \cite{GhMaAs09}. Moreover, it aims at getting a critical cost advantage, since the deployment of wireless transmitters is generally simpler and less expensive than that of optical fibers.

In this paper, we provide an original optimization model for the design of 3-architecture urban access networks integrating wireless and wired connections. A distinctive feature of our model w.r.t. state-of-the-art literature available on the topic (see \cite{LeEtAl13} and  \cite{GrRaWe14} for an overview) is to include the mathematical expressions that model wireless signal coverage and that evaluates the relation between useful and interfering signals.
%(the widely known, signal-to-interference ratios \cite{Ra01}).
The inclusion of such expressions is critical in any wireless network design problem considering wireless signal coverage, since the exclusion may lead to wrong design decisions (see \cite{DA10,DA15} for a discussion). The resulting problem has proved to be very difficult to solve even for a state-of-the-art commercial MIP solver like IBM ILOG CPLEX \cite{CPLEX}.

\smallskip
\noindent
In this work, our main original contributions are:
\begin{enumerate}
  \item we propose the first optimization model for the problem of optimally designing a 3-architecture access network, explicitly modelling the \emph{signal-to-interference formulas} that express wireless signal coverage. Specifically, we trace back the design problem to a 3-architecture variant of the \emph{Connected Facility Location Problem} that includes additional variables and constraints for modelling the service coverage of the wireless architecture;
  \item in order to strengthen the mathematical formulation of the problem, we propose to include two families of valid inequalities that model conflicts between variables representing the activation of wireless transmitters and the assignment of users to the transmitters;
  \item we develop a new primal heuristic for solving the problem. The heuristic is based on the combination of a probabilistic variable fixing procedure, guided by suitable \emph{Linear Programming} (LP) relaxations of the problem, with an exact \emph{Mixed Integer Programming} (MIP) heuristic, which provides for executing a \emph{very large neighborhood search} formulated as a \emph{Mixed Integer Linear Program} (MILP) and solved exactly by a state-of-the-art MIP solver;
  \item we present computational experiments over a set of realistic network instances, showing that our new algorithm is able to produce solutions of much higher quality than those returned by a state-of-the-art MIP solver.
\end{enumerate}
The remainder of this paper is organized as follows: in Section 2, we review the 2-architecture Connected Facility Location Problem; in Section 3, we introduce the new formulation for 3-architecture network design; in Sections 4 and 5, we present our new metaheuristic and discuss computational results.

\section{2-Architecture Connected Facility Location}

We start our modeling considerations by taking into account a generalization of a \emph{Connected Facility Location Problem} (ConFL) including two types of architectures. For an exhaustive introduction to foundations of network flow theory on graphs and to the ConFL, we refer the reader to \cite{AhMaOr93} and \cite{GoLj11}. The ConFL can be essentially described as the problem of a) deciding the assignment of a set of served users to a set of open facilities and b) how to connect open facilities through a Steiner tree, in order to minimize the total cost deriving from opening and connecting facilities and the assignment of facilities to users. The ConFL has been introduced and proven to be NP-Hard in \cite{GuEtAl01}. A hop-constrained version of ConFL that is related to the design of single-architecture access network has been studied in \cite{LiGo12}.

The canonical ConFL considers a \emph{single architecture} and can be associated to the design of an urban access network using a single \emph{technology} (i.e., either optical fiber or copper connections).
However, as we highlighted in the introduction, a new modern trend is to integrate two architectures and mix optical fiber and copper connection technologies. This leads to an extension of the ConFL that has been first considered and modeled in \cite{LeEtAl13} and that we denote by 2-ConFL. We now proceed to define an optimization model for 2-ConFL that we use as basis for introducing our new model for a 3-architecture ConFL including wireless technology.

The 2-ConFL in access network design involves a set of potential facilities that can install one among two technologies and that provide a telecommunication service to a set of potential users. A served user must be assigned to exactly one open facility and each open facility must be connected to a central office. The aim is to guarantee a minimum coverage of users by each technology, while minimizing the cost of deployment of the network.

To formally define the 2-ConFL, it is useful to consider a modeling of the network as a directed graph $G(V,A)$ where:
\begin{enumerate}
\item the set of nodes $V$ is the (disjoint) union of i) a set of users $U$ associated with a weight $w_u \geq 0$ representing the importance of each user $u \in U$, ii) a set of facilities $F$ with opening cost $c_f^{t} \geq 0$, $\forall f \in F$ that depends upon the technology $t \in T$ used by $f$, iii) a set of central offices $\Gamma$, with opening cost $c_\gamma \geq 0$, $\forall \gamma \in \Gamma$, iv) a set of Steiner nodes $S$. We call \emph{core nodes} the subset of nodes $V^{\text{C}} = F \cup \Gamma \cup S$ that does not include the user nodes. Additionally, we denote by $F_{u}^{t}$ the subset of facilities using technology $t$ that may serve user $u$ and by $U_f^{t}$ the subset of users that may be served by facility $f$ when using technology $t$.
\item the set of arcs $A$ is the (disjoint) union of i) a set of \emph{core arcs} $A^{\text{C}} = \{(i,j): i,j \in V^{\text{C}}\}$ that represent connections only between core nodes and are associated with a cost of realization $c_{ij} \geq 0$; iii) a set of \emph{assignment arcs} $A^{\text{ASS}} = \{(f,u) \in A: u \in U, f \in F_u\}$ representing connection of facilities to users and associated with a cost of realization $c_{fu}^{t}$ that depends upon the used technology.
\end{enumerate}

\noindent
We call \emph{core graph} the subgraph $G^{\text{C}} (V^{\text{C}}, A^{\text{C}})$ of $G(V,A)$ representing the potential topology of the optical fiber deployment (\emph{core network}) that has the \emph{core nodes} as set of nodes and the \emph{core arcs} as set of arcs.  To take into account the opening cost of central offices, we use the common trick to add an artificial root node $r$ to $G(V,A)$ that is connected to each central office $\gamma \in \Gamma$ by an arc $(r,\gamma)$ associated with cost $c_{r\gamma}$ that is set equal to the cost $c_{\gamma}$ of opening $\gamma$. This entails the inclusion in $G(V,A)$ of an additional set of (artificial) arcs $A^{\text{R}} = \{(r,\gamma): \gamma \in \Gamma\}$. In what follows, we will use the notation $A^{\text{R-C}} = A^{\text{R}} \cup A^{\text{C}}$ to denote the union of the root and the core arcs.
The total cost of deployment of the access network is obtained by summing the cost of opening central offices and facilities, the cost of connections established within the core graph and the cost of connecting open facilities to served users.

For each architecture, it is necessary to ensure a minimum weighted coverage of users. Given the total weight of users $W = \sum_{u \in U} w_u$, we express the coverage requirement for the architecture corresponding to technology $t \in T$ by introducing thresholds $W_t \in [0,W]$, $t \in T$. We assume that $W_1 \leq W_2$, i.e. the coverage requirement of the more performing and costly technology $t = 1$ is not higher than that of the lower class technology $t = 2$. We base this on the realistic assumption that just a part of the users is willing to pay more for getting a higher quality of service.
%, whereas m thus is satisfied with the less performing (but also cheaper) technology.

On the basis of the previous formalization of the problem, we can finally introduce a \emph{mixed integer linear program} to model the 2-ConFL. To this end, we introduce the following family of variables: 1) facility opening variables $z_{f}^{t} \in \{0,1\}$ $\forall f \in F, t \in T$ - the generic $z_{f}^{t}$ is equal to $1$ if facility $f$ is open and uses technology $t$ and is $0$ otherwise; 2) arc installation variables $x_{ij} \in \{0,1\}$ $\forall (i,j) \in A^{\text{R-C}}$ - the generic $x_{ij}$ is equal to $1$ if the root or core arc $(i,j)$ is installed and is $0$ otherwise; 3) assignment arc variables $y_{fu}^{t} \in \{0,1\}$
%$\forall (f,u) \in A^{\text{ASS}}$
$\forall u \in U$, $t \in T$, $f \in F_u^{t}$ - the generic $y_{fu}^{t}$ is equal to $1$ if facility $f$ is connected to user $u$ by technology $t$ and is $0$ otherwise; 4) user variables $v_{u}^{t} \in \{0,1\}$, $\forall u \in U$, $t \in T$ - the generic $v_{u}^{t}$ is equal to $1$ if user $u$ is served by technology $t$ and is $0$ otherwise; 5) flow variables $\phi_{ij}^{f}$, $\forall (i,j) \in A^{\text{R-C}}$, $f \in F$ representing the amount of flow sent on a root or core arc (i,j) for facility $f$. The Mixed Integer Linear Program for 2-ConFL (2-ConFL-MILP) is then:
{
%\small
\begin{align}
\min &
\sum_{(i,j) \in A^{\text{R-C}}} c_{ij} \hspace{0.05cm} x_{ij}
+ \sum_{f \in F} \sum_{t \in T} c_{f}^{t} \hspace{0.05cm} z_{f}^{t}
+ \sum_{u \in U} \sum_{t \in T} \sum_{f \in F_u^{t}} c_{fu}^{t} \hspace{0.05cm} y_{fu}^{t}
&&
\mbox{(2-ConFL-MILP)}
\nonumber
\\
&
\sum_{t \in T} z_{f}^{t} \leq 1
&&
f \in F
\label{2CFL_singleTechFacility}
\\
&
\sum_{f \in F_u^{t}} y_{fu}^{t} = \hspace{0.1cm} v_u^{t}
&&
u \in U, t \in T
\label{2CFL_user}
\\
&
y_{fu}^{t} \leq z_{f}^{t}
&&
u \in U, f \in F, t \in T
\label{2CFL_linkAssignmentUser}
\\
&
\sum_{u \in U} \sum_{\tau = 1}^{t} w_u \hspace{0.1cm} v_u^{t} \geq W_t
&&
t \in T
\label{2CFL_weight}
\\
&
\sum_{(j,i) \in A^{\text{R-C}}} \phi_{ji}^{f} - \sum_{(i,j) \in A^{\text{R-C}}} \phi_{ij}^{f}
=
\left\{
    \begin{array}{l}
        - \sum_{t \in T} z_{f}^{t}
        \\
        0
        \\
        + \sum_{t \in T} z_{f}^{t}
    \end{array}
    \hspace{0.1cm}
    \begin{array}{l}
        \mbox{if $i=r$}
        \\
        \text{if $i \neq r,f$}
        \\
        \text{if $i=f$}
    \end{array}
\right.
&&
i \in V^{\text{C}} \cup \{r\}, f \in F
\label{2CFL_flowConservation}
\\
&
0 \leq \phi_{ij}^{f} \leq x_{ij}
&&
(i,j) \in A^{\text{R-C}}, f \in F
\label{2CFL_flowBounds}
\\
&
v_u^{t}, \hspace{0.1cm} z_f^{t}, \hspace{0.1cm} x_{ij}^{t}, \hspace{0.1cm} y_{fu}^{t} \in \{0,1\}
&&
(i,j) \in A, u \in U, f \in F, t \in T
\nonumber
\end{align}
}

\noindent
The objective function aims at minimizing the total cost, expressed as the sum of the cost of activating root and core arcs (note that the corresponding summation includes the cost of activated central offices, opened facilities and of activated assignment arcs).
The constraints \eqref{2CFL_singleTechFacility} impose that each facility is activated on a single technology, whereas constraints \eqref{2CFL_user} impose that if a user $u$ is served by technology $t$, exactly one of the assignment arcs coming from a facility that can serve $u$ is activated on technology $t$. The constraints \eqref{2CFL_linkAssignmentUser} link the opening of a facility $f$ on technology $t$ to the activation of assignment arcs involving $f$ and $t$. The constraints \eqref{2CFL_weight} impose the coverage requirement for each user (note that here the weighted sum of users getting the better technology $t = 1$ contributes to satisfying the requirement for the coverage of the worse technology).

The constraints \eqref{2CFL_flowConservation} and \eqref{2CFL_flowBounds} jointly model the fiber connectivity within the core network as a multicommodity flow problem that includes one commodity per facility. Specifically, \eqref{2CFL_flowConservation} represents flow conservation in root and core nodes, while \eqref{2CFL_flowBounds} are variable upper bound constraints that express the linking between the activation of a root or core arc and the activation of the arc.

We note that in contrast to the formulation proposed in \cite{LeEtAl13}, which models connectivity within the core network by cut-set inequalities and whose size is thus potentially exponential in the size of the problem input, we adopt a \emph{compact formulation}  based on multicommodity flows that is polynomial in the size of the problem input. The compact formulation is indeed more suitable for being used in our new heuristic, not requiring the execution of additional time consuming separation routines.

\section{3-Architecture Connected Facility Location}

We now proceed to introduce our new original generalization of the 2-ConFL problem, which additionally considers wireless \emph{FTTA architecture} and explicitly embed the formulas expressing wireless coverage for a user.

As first step, we need to add an additional element to the set of available technologies, i.e. $T := T \cup \{3\}$ with index $t = 3$ denoting the wireless technology. We then assume that each facility $f \in F$ can also accommodate a \emph{wireless transmitter}, which may provide service connection \emph{without need of cables} to a subset of users. Transmitters are characterized by a number of radio-electrical parameters to set (e.g., the power emission, the frequency channel used to transmit, the modulation and coding scheme - see \cite{Ra01}). In principle, all these parameters can be set in an optimal way, by expressing their setting through a suitable mathematical optimization problem. However, just a (small) subset of parameters are typically optimized in a wireless network design problem \cite{DA10,DA11}. A decision that is included in practically every design problem is the setting of power emissions. This is indeed a crucial decision that deeply influences the possibility of covering users with service \cite{DAMaSa13}.

In order to model the power emission of a wireless facility $f \in F$, we introduce a semi-continuous power variable $p_f \in [P_{\min},P_{\max}] \hspace{0.2cm} \forall \hspace{0.1cm} f \in F$. A user $u$ receives power from each wireless facility $f \in F$ and the power $P_{f}(u)$ that $u$ receives from $f$ is proportional to the power emitted by $f$ by a factor $a_{fu} \in [0,1]$, i.e. $P_{f}(u) = a_{fu}\cdot p_f$. The factor $a_{fu}$ is called {\em fading coefficient} and
expresses the reduction in power that a signal
propagating from $f$ to $u$ experiences \cite{Ra01}.
We say that a user $u \in U$ is {\em covered} or {\em served} if it receives the wireless service signal within a minimum level of quality. The service is provided by one single transmitter, chosen as \emph{server} of the user, while all the other transmitters interfere with the server and reduce the quality of service.
The minimum quality condition can be expressed through the \emph{Signal-to-Interference Ratio (SIR)}, a measure comparing the power received from the server with the sum of the power received by the interfering transmitters \cite{Ra01}:
\begin{equation}
\label{eq:firstSIRineq} \frac{a_{fu} \hspace{0.1cm}
p_f}{\eta + \sum_{k \in F\setminus\{f\} } a_{ku}
\hspace{0.1cm} p_k}
\hspace{0.1cm}
\geq
\hspace{0.1cm}
\delta \;  .
\end{equation}

\noindent
The user is served if the SIR is above a threshold $\delta > 0$ that depends upon the wanted quality of service.
We remark that in the denominator we must also include a constant $\eta > 0$ representing the noise of the system. By simple linear algebra operations, inequality \eqref{eq:firstSIRineq} can be transformed in the so-called \emph{SIR inequality}:
%
%\begin{equation}\label{eq:secondSIRineq}
    $
    a_{fu} \hspace{0.1cm} p_f - \delta \sum_{k \in F\setminus\{f\} } a_{ku}
    \hspace{0.1cm} p_k
    \hspace{0.1cm}
    \geq
    \hspace{0.1cm}
    \delta \hspace{0.1cm} \eta \;  .
    $
%\end{equation}
Since we do not know in advance which wireless facility $f \in F$ will be the server of user $u \in U$ (establishing the assignment facility-user is part of the decision process), given a user $u \in U$ we have one SIR inequality
%(\ref{eq:secondSIRineq})
for each potential server $f \in F$, which must be activated or deactivated depending upon the assignment.
In order to ensure that $u$ is served through a wireless connection, at least one SIR inequality
%\eqref{eq:secondSIRineq}
must be satisfied. We are thus actually facing a disjunction of constraints, which, according to a standard approach of Mixed Integer Programming (see \cite{NeWo88}), can be represented by a variant of the SIR inequality
%\eqref{eq:secondSIRineq}
that includes a sufficiently large positive constant $M$ (the so-called \emph{big-M coefficient}) and the assignment variable $y_{fu}^{3}$ representing the service connection of $u$ through facility $f$ by technology $t = 3$, namely:
\begin{equation}\label{eq:SIR-BIGM}
a_{f u} p_{f} - \delta \sum_{k \in F \backslash \{f\}} a_{ku} p_{k}
+ M (1 - y_{fu}^{3})
\geq \delta N
\end{equation}

\noindent
It is immediate to check that if
$y_{fu}^{3} = 1$, then $u$ is wirelessly served by $f$ and (\ref{eq:SIR-BIGM}) reduces to a SIR inequality to be satisfied. If instead $y_{fu}^{3} = 0$, then $u$ is not wirelessly served by $f$ and $M$ activates, thus making (\ref{eq:SIR-BIGM}) redundant and satisfied by any power vector $(p_{1}, p_{2}, \ldots, p_{|F|})$.
The MILP for 3-ConFL is then:
{
%\small
\begin{align}
\min &
\sum_{(i,j) \in A^{\text{R-C}}} c_{ij} \hspace{0.05cm} x_{ij}
+ \sum_{f \in F} \sum_{t \in T} c_{f}^{t} \hspace{0.05cm} z_{f}^{t}
+ \sum_{u \in U} \sum_{t \in T} \sum_{f \in F_u^{t}} c_{fu}^{t} \hspace{0.05cm} y_{fu}^{t}
&&
\mbox{(3-ConFL-MILP)}
\nonumber
\\
&
\sum_{t \in T} z_{f}^{t} \leq 1
&&
f \in F
\nonumber
%\label{3CFL_singleTechFacility}
\\
&
\sum_{f \in F_u^{t}} y_{fu}^{t} = \hspace{0.1cm} v_u^{t}
&&
u \in U, t \in T
\nonumber
%\label{3CFL_user}
\\
&
y_{fu}^{t} \leq z_{f}^{t}
&&
u \in U, f \in F, t \in T
\nonumber
%\label{3CFL_linkAssignmentUser}
\\
&
\sum_{u \in U} \sum_{\tau = 1}^{t} w_u \hspace{0.1cm} v_u^{t} \geq W_t
&&
t \in T
\nonumber
%\label{3CFL_weight}
\\
&
\sum_{(j,i) \in A^{\text{R-C}}} \phi_{ji}^{f} - \sum_{(i,j) \in A^{\text{R-C}}} \phi_{ij}^{f}
=
\left\{
    \begin{array}{l}
        - \sum_{t \in T} z_{f}^{t}
        \\
        0
        \\
        + \sum_{t \in T} z_{f}^{t}
    \end{array}
    \hspace{0.1cm}
    \begin{array}{l}
        \mbox{if $i=r$}
        \\
        \text{if $i \neq r,f$}
        \\
        \text{if $i=f$}
    \end{array}
\right.
&&
i \in V^{\text{C}} \cup \{r\}, f \in F
\nonumber
%\label{3CFL_flowConservation}
\\
&
0 \leq \phi_{ij}^{f} \leq x_{ij}
&&
(i,j) \in A^{\text{R-C}}, f \in F
\nonumber
%\label{3CFL_flowBounds}
\\
&
a_{f u} p_{f} - \delta \sum_{k \in F \backslash \{f\}} a_{ku} p_{k}
+ M (1 - y_{fu}^{3})
\geq \delta N
&&
f \in F, u \in U
\label{3CFL_SIR constraints}
\\
&
0 \leq P^{\min} z_{f}^{3} \leq p_f \leq P^{\max} z_{f}^{3}
&&
f \in F
\label{3CFL_powerVariables}
\\
&
v_u^{t}, \hspace{0.1cm} z_f^{t}, \hspace{0.1cm} x_{ij}^{t}, \hspace{0.1cm} y_{fu}^{t} \in \{0,1\}
&&
(i,j) \in A, u \in U, f \in F, t \in T
\nonumber
\end{align}
}

\noindent
The major modifications w.r.t. the formulation (2-ConFL-MILP) concern: 1) the introduction of the variable bound constraints \eqref{3CFL_powerVariables} that express the semi-continuous nature of variables $p_f$ (when $z_{f}^{3}=0$, facility $f$ does not install a wireless transmitter and the power $p_f$ is thus forced to 0; when instead $z_{f}^{3}=1$, the transmitter is installed and its power must lie in $[P^{\min},P^{\max}]$); 2) the introduction of the SIR constraints \eqref{3CFL_SIR constraints} for expressing the wireless coverage conditions.

\subsubsection{Strengthening 3-ConFL-MILP.}
\label{subsec:ACO}

A key ingredient of the probabilistic fixing that we adopt in our new heuristic is represented by the combination of an a-priori and an a-posteriori measure of fixing attractiveness based on linear relaxations of 3-ConFL-MILP. In particular, we obtain the a-priori measure considering a \emph{tighter formulation} (informally speaking, a problem with a ``mathematically stronger'' structure) defined by adding two class of valid inequalities to 3-ConFL-MILP: 1) superinterferer inequalities; 2) conflict inequalities. These two families of inequalities were respectively introduced in \cite{DA10} and \cite{DAMaSa11} and we refer the reader to these works for a detailed description. Here, we just provide a concise introduction to them.

The first class of inequalities captures the existence of so-called \emph{superinterferers}: a superinterferer is an interfering transmitter that alone can deny service coverage to a user even when it emits at minimum power and the serving transmitter emits at maximum power. The corresponding valid inequalities are logical constraints of the form:
\begin{equation}
y_{fu}^{3} \leq 1 - z_k^{3}
\hspace{0.5cm}
\forall \hspace{0.1cm} k \in K \backslash \{f\}: k \text{ is superinterfer for $u$ served by $f$}
\label{eq:superinterfer}
\end{equation}

\noindent
expressing that if $k$ is a superinterferer facility and is activated, then the variable assigning user $u$ to $f$ installing wireless technology is forced to $0$, since the corresponding SIR constraint cannot be satisfied (notice that the set of superinterferers depends upon the considered user and the user serving facility).

The second class of valid inequalities captures the existence of couples of SIR constraints that involve just two wireless facilities and that cannot be satisfied at the same time.
More formally, consider the two SIR constraints corresponding to two users $u_1, u_2$ served by two distinct wireless facilities $f_1, f_2$, namely:
1) $a_{f_1 u_1} p_{f_1} - \delta a_{f_2 u_1} p_{f_2} \geq \delta N;$
2) $a_{f_1 u_1} p_{f_1} - \delta a_{f_2 u_1} p_{f_2} \geq \delta N$
(respectively representing $u_1$ served by $f_1$ and interfered by $f_2$ and $u_2$ served by $f_2$ and interfered by $f_1$).
If there is no power vector $(p_1, p_2)$ that satisfies the power bounds \eqref {3CFL_powerVariables} and the two SIR constraints,
%\begin{eqnarray*}
%a_{f_1 u_1} p_{f_1} - \delta a_{f_2 u_1} p_{f_2} \geq \delta N
%\\
%a_{f_2 u_2} p_{f_2} - \delta a_{f_1 u_2} p_{f_1} \geq \delta N
%\end{eqnarray*}
then the following is a \emph{valid inequality} for 3-ConFL-MILP stating that both SIR constraints cannot be activated simultaneously:
\begin{equation}
y_{f_1 u_1}^{3} + y_{f_2 u_2}^{3} \leq 1
\label{eq:conflict}
\end{equation}

\noindent
Such valid inequalities can be easily identified in a pre-processing phase and can be added to the formulation to get remarkable strengthening (see \cite{DAMaSa11}).
In the next section, we denote by \emph{Strong-3-ConFL-MILP}, the problem 3-ConFL-MILP suitably strengthened by inequalities \eqref{eq:superinterfer} and \eqref{eq:conflict}.

\section{A primal heuristic for the 3-ConFL-MILP}
\label{sec:ACO}

Being a mixed integer linear program, the problem 3-ConFL-MILP could in principle be solved by using a commercial MIP solver, such as IBM ILOG CPLEX \cite{CPLEX}. However, the introduction of the wireless technology and of constraints \eqref{3CFL_SIR constraints} make 3-ConFL-MILP a very challenging extension of the 2-ConFL problem that result very difficult even for state-of-the-art solvers. According to our direct experience on realistic instances, in many cases CPLEX had big difficulties in finding good quality solutions even after several hours of computations. We observed analogue computational difficulties also in other problems based on the combination of flow models with signal-to-interference constraints (e.g., \cite{DeDAKa15}).

As an alternative to the direct use of a MIP commercial solver, we thus propose a new heuristic that combines a \emph{probabilistic fixing procedure}, guided by the solution of peculiar linear relaxations of 3-ConFL-MILP, with an MIP heuristic, based on an \emph{exact very large neighborhood search}.
The probabilistic fixing is partially inspired by the algorithm ANTS (\emph{Approximate Nondeterministic Tree Search}) \cite{Ma99} a refined version of an ant colony algorithm that tries to exploit information about bounds available for the optimization problem. More precisely, our new heuristic is based on considerations about the use of linear relaxations in place of generic bounds that have been first made in \cite{DAKrPu15} and \cite{DANa15}.

Since we exploit linear relaxations, in contrast to ``simple'' heuristics, we can provide a certificate of quality for the best solution produced by our heuristic. The certificate assumes the form of an \emph{optimality gap}, measuring how far the best solution is from the best lower bound given by Strong-3-ConFL-MILP.

It is nowadays widely known that Ant Colony Optimization (ACO) is a metaheuristic inspired by the behaviour of ants looking for food, initially proposed by Dorigo and colleagues for combinatorial optimization and then extended and refined in many works (e.g., \cite{DoMaCo96,GaMoWe12,Ma99,BlEtAl11}% - see also \cite{BlEtAl11} for an overview).
ACO is essentially centered on the execution of a cycle where a number of feasible solutions are iteratively built, using information about solution construction executed in previous runs of the cycle. An ACO algorithm (ACO-alg) presents the general structure of Algorithm \ref{(Gen-ACO)}.
\begin{algorithm}
\caption{General ACO Algorithm (ACO-alg)}
\label{(Gen-ACO)}
\begin{algorithmic}[1]
\While{an arrest condition is \emph{not} satisfied}
    \State ant-based solution construction
    \State pheromone trail update
\EndWhile
\State local search
\end{algorithmic}
\end{algorithm}

In the step 2 of the while-cycle, a number of \emph{ants} are defined and each ant builds a feasible solution in an iterative way. At every iteration, the ant is in a \emph{state} that corresponds with a \emph{partial solution} and can further complete the solution by making a \emph{move}. The move corresponds to fixing the value of a not-yet-fixed variable and is chosen in a probabilistic way, evaluating a measure that combines an \emph{a-priori} and an \emph{a-posteriori} measure of fixing attractiveness. The a-priori attractiveness measure is called \emph{pheromone trail value} in an ACO-alg context and is updated at the end of the construction phase, in the attempt of rewarding good fixing and penalizing bad fixing. Once that an  arrest condition is met (typically, reaching a time limit), a local search is executed to improve the quality of the produced solutions and possibly identify a local optimum.

We stress that the algorithm that we propose is \emph{not} an ACO-alg, but is rather an evolution and refinement of the ANTS algorithms that we strengthen by the use of peculiar linear relaxations. Specifically, in our case, the a-priori measure is provided by a strengthened linear relaxation of the problem - we use Strong-3-ConFL-MILP, namely problem 3-ConFL-MILP strengthened by inequalities \eqref{eq:superinterfer} and \eqref{eq:conflict} - whereas the a-posteriori measure is provided by the linear relaxation of 3-ConFL-MILP for partial fixing of the facility opening variables.
\\
We now proceed to describe in detail our new primal heuristic.

\medskip
\noindent
\textbf{Feasible solution construction.}
%\subsubsection{Feasible solution construction.}
%\label{subsubsec:solConstruction}
To explain how we build a feasible solution for the 3-ConFL-MILP, we first introduce the concept of \emph{Facility Opening state}:
\begin{definition}
  Facility opening state (FOS): let $F \times T$ be the set of couples $(f,t)$ that represent the activation of a facility $f$ on a technology $t$. An FOS specifies an opening of a subset of facilities $\bar{F} \subseteq F$ on some technologies and excludes that the same facility is opened on more than one technology (i.e., $FOS \subseteq F \times T:
   $\hspace{0.1cm}$
   \not\exists (f_1,t_1), (f_2,t_2) \in FOS:
   \hspace{0.1cm}
   f_1 = f_2$ and $t_1 \neq t_2$).
\end{definition}

\noindent
Given a FOS and a facility-technology couple $(f,t) \in FOS$, we denote by
$W^{\text{POT}}_{ft}$ the total weight of users that can be potentially served by $f$ activated on technology $t$, i.e. $W^{\text{POT}}_{ft} = \sum_{u \in U_{f}^{t}} w_u$.
We introduce this measure to distinguish between a partial and complete FOS for a technology $t \in T$. We say that a FOS is \emph{partial for technology $t$} when the total weight of potential users that can be served by facilities appearing in the $FOS$ using technology $t$ does not reach the minimum coverage requirements $W_t$ for $t$, i.e.:
\begin{eqnarray}
\label{eqn:completeness}
\sum_{f \in F: (f,t) \in FOS} \hspace{0.1cm} \sum_{u \in U_{f}^{t}} w_u < W_t \; .
\end{eqnarray}

\noindent
On the contrary, we say that a FOS is \emph{complete for technology $t$} when the total weight is not lower than $W_t$.
Additionally, we call \emph{fully complete} a FOS that is complete for all technologies $t \in T$.
We introduce the concept of completeness and the formula \eqref{eqn:completeness} in order to guide and limit the probabilistic fixing of facility opening variables during the construction phase of feasible solutions.

Given a \emph{partial} FOS for technology $t$, the probability $p_{ft}^{\text{FOS}}$ of operating an additional fixing $(f,t) \not\in FOS$, thus making a further step towards reaching a complete FOS, is set according to the formula:
\begin{equation}
\label{probability}
p_{ft}^{\text{FOS}} = \frac{\alpha \hspace{0.1cm} \tau_{ft} + (1-\alpha) \hspace{0.1cm} \eta_{ft}}
                {\sum_{(k,t) \not\in FOS} \alpha \hspace{0.1cm} \tau_{kt} + (1-\alpha) \hspace{0.1cm} \eta_{kt}} \; ,
\end{equation}

\noindent
which provides for a convex combination of the a-priori attractiveness measure $\tau_{ft}$ and the a-posteriori attractiveness measure $\eta_{ft}$ through factor $\alpha \in [0,1]$. In our specific case, $\tau_{ft}$ is provided by the optimal value of the linear relaxation Strong-3-ConFL-MILP including the additional fixing $z_{f}^{t} = 1$, whereas $\eta_{kt}$ is the value of the linear relaxation of 3-ConFL-MILP obtained for a specified partial fixing of the facility opening variables $z$.
We remark that \eqref{probability} is a revised formula that was proposed in \cite{Ma99} to improve the computationally inefficient canonical formula of ACO, which includes products and powers of measures and depends upon a higher number of parameters.

At the end of a solution construction phase, the a-priori measures $\tau$ are updated, evaluating how good the fixing resulted in the obtained solutions. We stress that for the update we do not rely on the canonical ACO formula including the pheromone evaporation parameter, whose setting may result very tricky, but we use a revised version of the improved formula proposed for ANTS in \cite{Ma99}.
To define the new formula, we first introduce the concept of \emph{optimality gap} (\emph{OGap}) for a feasible solution of value $v$ and a lower bound $L$ that is available on the optimal value $v^{*}$ of the problem (note that it holds $L \leq v^{*} \leq v$): the \emph{OGap} provides a measure of the quality of the feasible solution, comparing its value to the lower bound and is formally defined as $OGap(v,L)$ $= (v - L)/v$. The a-priori attractiveness measure that we use is:
\begin{equation}
\small
\label{updateFormula}
\tau_{ft}(h) = \tau_{ft}(h-1)
\hspace{0.05cm} + \hspace{0.05cm} \sum_{\sigma=1}^{\Sigma} \Delta \tau_{ft}^\sigma
%\hspace{0.2cm}
\mbox{ with }
\Delta \tau_{ft}^\sigma =
\tau_{ft}(0)
\cdot \left(
\frac{OGap(\bar{v},L) - OGap(v_{\sigma},L)}{OGap(\bar{v},L)}
\right)
%\; ,
\end{equation}

\noindent
where $\tau_{ft}(h)$ is the a-priori attractiveness of fixing $(f,t)$ at fixing iteration $h$, $L$ is a lower bound on the optimal value of the problem (we remember that as lower bound we use the optimal value of the strengthened formulation Strong-3-ConFL-MILP),
$v_{\sigma}$ is the value of the $\sigma$-th feasible solution built in the last construction cycle and $\bar{v}$ is the (moving) average of the values of the $\Sigma$ solutions produced in the previous construction phase. $\Delta \tau_{ft}^\sigma$ represents the penalization/reward factor for a fixing and depends upon the initialization value $\tau_{ft}(0)$ of $\tau$ (in our case, based upon the linear relaxation of Strong-3-ConFL-MILP), combined with the relative variation in the optimality gap that $v_{\sigma}$ implies w.r.t. $\bar{v}$.
We note that the use of a relative gap difference in \eqref{updateFormula} allows us to reward or penalize fixing adopted in the last solution making a comparison with the average quality of the last $\Sigma$ solutions constructed.

Once that a \emph{fully complete FOS} is built, we have characterized an opening of facilities that can \emph{potentially} satisfy the requirements on the weighted coverage for each technology. We use the term ``potentially'', since the activation of facilities specified by the FOS does not necessarily admit a feasible completion in terms of connectivity variables and assignment of users of facilities: it is indeed likely that not all the SIR constraints \eqref{3CFL_SIR constraints} corresponding to wireless facilities can be activated simultaneously because of interference effects.
It is thus possible that a complete FOS will result infeasible. Since a risk of infeasibility is present, after the construction of a complete FOS, we execute a \emph{check-and-repair phase}, in which the feasibility of the FOS is checked and, if not verified, we make an attempt to repair and make it feasible. The reparation attempt is based on the same MIP heuristic based on an exact very large neighborhood search that we adopt at the end of the construction phase to possibly improve a feasible solution (see the next subsection for details).

Given a FOS that is complete for all technologies, we check its feasibility and attempt at finding a feasible solution for the complete problem 3-ConFL-MILP by defining a restricted version of 3-ConFL-MILP, where we set $z_{f}^{t} = 1$ if $(f,t) \in FOS$. We solve this restricted problem through the MIP solver with a time limit: if this problem is recognized as infeasible by the solver, we run the MIP heuristic for reparation. Otherwise, we run the solver to possibly find a solution that is better than the best incumbent solution.

\medskip
\noindent
\textbf{MIP-VLNS - an exact MIP repair/improvement heuristic.}
%\subsubsection{MIP-VLNS - an exact MIP repair/improvement heuristic.}
%\label{subsubsec:MIP-VNS}
To repair an infeasible partial fixing of the variables $z$ induced by a complete FOS or to improve an incumbent feasible solution, we rely on an MIP heuristic that operates a very large neighborhood search \emph{exactly}, by formulating the search as a mixed integer linear program solved through an MIP solver \cite{BlEtAl11}.
Specifically, given a (feasible or infeasible) and possibly not complete fixing $\bar{z}$ of variables, we define the neighborhood ${\cal N}$ including all the feasible solutions of 3-ConFL-MILP that can be obtained by modifying at most $n > 0$ components of $\bar{z}$ and leaving the remaining variables free to vary. This condition can be expressed in 3-ConFL-MILP by adding an \emph{hamming distance constraint} imposing an upper limit $n$ on the number of variables in $z$ that change their value w.r.t. $\bar{z}$:
$$
%HD (\bar{z},z)
%\hspace{0.1cm} = \hspace{0.1cm}
\sum_{(f,t) \in F \times T: \; \bar{z}_{f}^{t} = 0} z_{f}^{t} + \sum_{(f,t) \in F \times T: \; \bar{z}_{f}^{t} = 1} (1 - z_{f}^{t}) \hspace{0.1cm} \leq \hspace{0.1cm}
n
$$

\noindent
The modified problem is then solved through an MIP solver like CPLEX, running with a time limit. Imposing a time limit is essential from a practical point of view: optimally solving the exact search can take a very high amount of time to close the optimality gap; additionally, a state-of-the-art MIP solver is usually able to quickly find solutions of good quality for large problems whose size has been conveniently reduced by fixing. In what follows, we denote the overall procedure for repair/improvement that we have discussed by MIP-VLNS.

\medskip
\noindent
\textbf{The complete algorithm.}
%\subsubsection{The complete algorithm.}
%\label{subsubsec:algorithm}
The complete algorithm for solving the 3-ConFL-MILP is presented in Algorithm \ref{ALG_3-ConFL-MILP}. We base the algorithm on the execution of two nested loops: the outer loop runs until reaching a global time limit and contains an inner loop inside which we define $\Sigma$ feasible solutions, by first defining a complete FOS and then executing the MIP heuristic to repair or complete the fixing associated with the FOS.
More in detail, the first algorithmic task is to solve the linear relaxation of Strong-3-ConFL-MILP for each fixing $z_{f}^{t} = 1$, getting the corresponding optimal value and using it to initialize the a-priori measure of attractiveness $\tau_{ft}(0)$. This is followed by the definition of a solution $X^{*}$ that represents the best solution found during the execution of the algorithm. Each run of the inner loop provides for building a complete FOS by considering, in order, fiber, copper and wireless technology. The complete FOS is built according to the procedure using the probability measures \eqref{probability} and update formulas \eqref{updateFormula} that we have discussed before.
The complete FOS provides a (partial) fixing of the facility opening variables $\bar{z}$ and the MIP solver uses it as a basis for finding a complete feasible solution $X^{*}$ to the problem. If $\bar{z}$ is recognized as an infeasible fixing by the MIP solver, then we run the MIP-VLNS in a reparation mode. If instead $\bar{z}$ is feasible and leads to find a feasible solution to 3-ConFL-MILP that is better than the best solution found $X^{B}$ in the current run of the inner loop, then $X^{B}$ is updated. Then the inner loop is iterated. After that the execution of the inner loop is concluded, the a-priori measures $\tau$ are updated according to formula \eqref{updateFormula}, considering the quality of the produced solutions, and we check the necessity of updating the global best solution $X^{*}$. After having reached the global time limit, the heuristic MIP-VLNS is eventually run with the aim of improving the best solution found $X^{*}$.
\begin{algorithm}
\caption{
- Heuristic for 3-ConFL-MILP
}
\label{ALG_3-ConFL-MILP}
\begin{algorithmic}[1]
    \State compute the linear relaxation of Strong-3-ConFL-MILP for all $z_{f}^{t} = 1$ and initialize the values $\tau_{ft}(0)$ with the corresponding optimal values
    \State let $X^{*}$ be the best feasible solution found
    \While{a global time limit is not reached}
        \State let $X^{B}$ be the best solution found in the inner loop
        \For{$\sigma := 1$ to $\Sigma$}
            \State build a complete FOS
            \State solve 3-ConFL-MILP imposing the fixing $\bar{z}$ specified by the FOS
            \If {3-ConFL-MILP with fixing $\bar{z}$ is infeasible}
                \State run MIP-VLNS for repairing the fixing $\bar{z}$
            \EndIf
            \If {a feasible solution $\bar{X}$ is found by the MIP solver and $c(\bar{X}) < c(X^{B})$}
                \State update the best solution found $X^{B} := \bar{X}$
            \EndIf
        \EndFor
        \State update $\tau$ according to (\ref{updateFormula})
        %\State run MIP-VLNS($\bar{x},\bar{y}$) for improving $(x^{B},y^{B})$
        %
        \If {$c(X^{B}) < c(X^{*})$}
            \State update the best solution found $X^{*}:= X^{B}$
        \EndIf
    \EndWhile
    \State run MIP-VLNS for improving $X^{*}$
    \State return $X^{*}$
\end{algorithmic}
\end{algorithm}

\section{Computational results}
\label{sec:computations}

We tested the performance of our algorithm on 15 instances based on realistic network data defined within past consulting and industrial projects for a major telecommunication company. The experiments were performed on a 2.70 GHz Windows machine with 8 GB of RAM and using IBM ILOG CPLEX 12.5 as MIP solver. The code was written in C/C++ and is interfaced with CPLEX through Concert Technology. The experiments ran with a time limit of 3600 seconds. All the instances refers to a urban district in the metropolitan area of Rome (Italy) and considers different traffic generation and user location scenarios. The considered area has been discretized into a grid of about 450 pixels, following the \emph{testpoint model} recommended by international telecommunications regulatory bodies for  wireless signal evaluation (see \cite{AGCOM,DA10}). We considered 30 potential facility locations that can accommodate any of the 3 technology considered in the study and can be connected to 5 potential central offices. On the basis of past experience and preliminary tests, we imposed the following setting of the parameters of the heuristic: $\alpha = 0.5$ (a-priori and a-posteriori attractiveness are balanced), $\Sigma = 5$ (number of solutions built in the inner loop before updating the a-priori measure and  width of the moving average). Additionally, we imposed a time limit of 3000 seconds to the execution of the outer loop of Algorithm \ref{ALG_3-ConFL-MILP} and a limit of 600 seconds to the execution of the improvement heuristic MIP-VLNS.

The computational results are presented in Table \ref{table:results}: here, for each instance instance, we report its ID, the best percentage optimality gap \emph{Gap-CPLEX\%} reached by CPLEX within the time limit, the best percentage optimality gap reached by our heuristic within the time limit \emph{Gap-Heu\%}. In the case of the heuristic, we note that the gap is obtained combining the best feasible solution found by Algorithm \ref{ALG_3-ConFL-MILP} with the best known lower bound obtained by CPLEX using the strengthened formulation Strong-3-ConFL-MILP.
\begin{table}
\caption{Experimental results}
\label{table:results}
\begin{center}
\begin{tabular}{c | c c c}
\hline
\hspace{0.1cm} ID \hspace{0.1cm}
& \hspace{0.1cm} Gap-CPLEX\% \hspace{0.1cm}
& \hspace{0.1cm} Gap-Heu\% \hspace{0.1cm}
& \hspace{0.1cm} $\Delta$Gap\% \hspace{0.1cm}
\\ [1pt]
\hline
I1  &   148.57  &   131.23  &   -11.67\\
I2 	&   136.74  &   106.16  &   -22.36\\
I3 	&  	99.46   &   72.96   &   -26.64\\
I4 	&  	156.47	&   123.73  &   -20.92\\
I5 	&  	78.86	&   49.98   &   -36.62\\
I6 	&  	93.42	&   64.04   &   -31.44\\
I7 	&  	117.00 	&   82.05	&   -29.48\\
I8 	&  	 95.21	&   59.73   &   -37.26\\
I9 	&  	178.94	&  119.62   &   -33.15\\
I10 &  	 98.80	&   77.66   &   -21.39\\
I11 &  	 89.13	&   66.17   &   -25.76\\
I12 &   104.11	&   71.23   &   -31.58\\
I13 &  	 95.20	&   52.08   &   -45.29\\
I14 &  	112.44	&   82.48   &   -26.64\\
I15 &  	103.00	&   74.30   &   -27.86\\
\hline
\end{tabular}
\end{center}
\end{table}

Concerning the results, the first critical observation to be made is that 3-ConFL-MILP results very challenging even for a modern state-of-the-art solver like CPLEX: the minimum gap obtained for the majority of instances results far beyond 90\%. We believe that such difficulty in solving the problem is particularly due to the presence of the SIR constraints \eqref{3CFL_SIR constraints} associated with wireless coverage: pure wireless coverage problem constitutes indeed already very challenging optimization problems, as discussed in \cite{DA10}. In comparison to CPLEX, our heuristic is able to get always (much) better optimality gaps, that on average are 28\% lower than those produced by CPLEX and can reach even reductions over 35\% (we remind that decreasing the optimality gap is crucial to ``move towards'' identifying the optimal value of an optimization problem, see \cite{NeWo88}). We believe that this is a very promising performance and that the heuristic deserves further investigations to be enhanced.

\section{Conclusion and future work}  \label{sec:end}
We considered the design of 3-architecture urban access networks, combining the use of wired optical fiber- and copper-based connections with wireless connections. In literature, it has been suggested that this problem can be modeled as a simple generalization of  2-architecture Connected Facility Location Problems, by including an additional technology index. However, this is a simplistic generalization that neglects the interaction between signals emitted by distinct wireless transmitters
and that may possibly lead to service coverage plans not implementable in practice. As a remedy, we have proposed a new optimization model that
%, besides being based on a  3-architecture Connected Facility Location model,
also includes the variables and constraints modeling the power emissions of wireless transmitters and the signal-to-interference formulas that are recommended for evaluating wireless service coverage.
The resulting mixed integer linear program results very challenging even for a state-of-the-art commercial MIP solver like CPLEX, so
we have proposed a new heuristic based on the combination of a probabilistic variable fixing procedure, guided by suitable linear relaxations of the problem, and an exact very large neighborhood search. Computational experiments on a set of realistic instances indicate that our heuristic can provide solutions associated with much lower optimality gaps than those returned by CPLEX.
As future work, we plan to refine the construction phase of the heuristic, studying further formulation strengthening, and to integrate it with a branch-and-cut algorithm to improve the overall computational performance. Additionally, we will consider data-uncertain of the problem, using Multiband Robust Optimization \cite{BuDA12a,DA15}.

%\subsection*{Acknowledgment}
%%
%\noindent
%{\small
%The work of Fabio D'Andreagiovanni and Jonad Pulaj was partially supported by the \emph{Einstein Center for Mathematics Berlin} (ECMath) through Project MI4 (ROUAN) and by the \emph{German Federal Ministry of Education and Research} (BMBF) through Project VINO (Grant 05M13ZAC) and Project \emph{ROBUKOM} (Grant 05M10ZAA)\cite{BaEtAl14}.
%}

\bibliographystyle{splncs}
\bibliography{DAndreagiovanni}

\end{document}